\documentclass{amsart}
\usepackage{amsmath}%
\usepackage{amsfonts}%
\usepackage{amssymb}%
\usepackage{graphicx}

\usepackage{hyperref}
\usepackage{enumerate}
\newtheorem{theorem}{Theorem}
\theoremstyle{plain}

\newtheorem{claim}{Claim}

\newtheorem{conjecture}{Conjecture}
\newtheorem{corollary}{Corollary}

\newtheorem{definition}{Definition}
\newtheorem{example}{Example}

\newtheorem{lemma}{Lemma}

\newtheorem{proposition}{Proposition}
\newtheorem{remark}{Remark}

\numberwithin{equation}{section}

\newcommand{\R}{\mathbb{R}}
\newcommand{\Sp}{\mathbb{S}}
\newcommand{\Q}{\mathbb{Q}}

\newcommand{\virg}{\mbox{, }}
\newcommand{\al}{\alpha}

\newcommand{\Ric}{\mbox{\normalfont{Ric}}}

\newcommand{\beeq}{\begin{equation}}
\newcommand{\eneq}{\end{equation}}
\newcommand{\beeqs}{\begin{eqnarray*}}
\newcommand{\eneqs}{\end{eqnarray*}}
\newcommand{\besp}{\begin{split}}
\newcommand{\ensp}{\end{split}}
\newcommand{\bepr}{\begin{proof}}
\newcommand{\enpr}{\end{proof}}
\newcommand{\beths}{\begin{theorem*}}
\newcommand{\enths}{\end{theorem*}}
\newcommand{\becor}{\begin{corollary}}
\newcommand{\encor}{\end{corollary}}
\newcommand{\bere}{\begin{remark}}
\newcommand{\enre}{\end{remark}}
\newcommand{\beres}{\begin{remark*}}
\newcommand{\enres}{\end{remark*}}
\newcommand{\bele}{\begin{lemma}}
\newcommand{\enle}{\end{lemma}}
\newcommand{\beles}{\begin{lemma*}}
\newcommand{\enles}{\end{lemma*}}
\newcommand{\bepro}{\begin{proposition}}
\newcommand{\enpro}{\end{proposition}}
\newcommand{\bepros}{\begin{proposition*}}
\newcommand{\enpros}{\end{proposition*}}
\newcommand{\becl}{\begin{claim}}
\newcommand{\encl}{\end{claim}}
\newcommand{\beex}{\begin{example}}
\newcommand{\enex}{\end{example}}
\newcommand{\beexs}{\begin{example*}}
\newcommand{\enexs}{\end{example*}}
\newcommand{\beco}{\begin{conjecture}}
\newcommand{\enco}{\end{conjecture}}
\newcommand{\becos}{\begin{conjecture*}}
\newcommand{\encos}{\end{conjecture*}}
\newcommand{\bede}{\begin{definition}}
\newcommand{\bedes}{\begin{definition*}}
\newcommand{\inj}{\mbox{\normalfont{inj}}}
\newcommand{\rad}{\mbox{\normalfont{\scriptsize{rad}}}}
\newcommand{\If}{\mbox{ if }}
\newcommand{\hess}{\mbox{\normalfont{Hess}}}
\newcommand{\N}{\mathbb{N}}
\newcommand{\grad}{\mbox{\normalfont{grad}}}
\newcommand{\Aand}{\mbox{ and }}
\textheight=\dimexpr
  \ifcase 10pt
    \or % 0
    \or % 1
    \or % 2
    \or % 3
    \or % 4
    \or % 5
    \or % 6
    \or % 7
    57\or % 8
    52\or % 9
    48\or % 10
    44\or % 11
    41\fi % 12
  \baselineskip+\topskip\relax
\begin{document}
\title[Submanifolds With Nonpositive Extrinsic Curvature]{Submanifolds With Nonpositive Extrinsic Curvature}

\author[S. Canevari]{Samuel Canevari}
\address{Samuel Canevari -- Universidade Federal de Sergipe \newline%
\indent Av. Ver. Ol\'{i}mpio Grande, Centro, 49500-000, Itabaiana, Brazil}
\email{samuel@ufs.br}%
%\urladdr{http://www.dmai.ufs.br/pagina/prof-samuel-cruz-canevari-5912.html}

\author[G. Freitas]{Guilherme Machado de Freitas}
\address{Guilherme M. de Freitas -- Instituto de Matem\'atica
Pura e Aplicada\newline%
\indent Estrada Dona Castorina, 110, Jardim Bot\^{a}nico, 22460-320,
Rio de Janeiro, Brazil}%
\email{gfreitas@impa.br}%
%\urladdr{https://www.facebook.com/guilherme.machadodefreitas}

\author[F. Manfio]{Fernando Manfio}
\address{Fernando Manfio -- Universidade de S\~ao Paulo \newline%
\indent Av. Trabalhador S\~{a}o-carlense, 400, Centro, 13560-970, 
S\~{a}o Carlos, Brazil}%
\email{manfio@icmc.usp.br}%
%\urladdr{ http://www.icmc.usp.br/pessoas/manfio}

\thanks{The first author's research was partially supported by CAPES/Brazil}
\thanks{The second author's research was partially supported by CNPq/Brazil}
\thanks{The third author's research was partially supported by FAPESP/Brazil, 
grant 2014/01989-9}

%\thanks{This paper is in final form and no version of it will be submitted for
%publication elsewhere.}
%\date{\today}
\subjclass[2010]{Primary 53C40, 53C42; Secondary 53A07, 53A35} %
\keywords{nonpositive extrinsic curvature, cylindrically bounded submanifolds, Otsuki's Lemma, Omori-Yau maximum principle}%
%\dedicatory{Dedicated to Professor X on the occasion of his seventieth birthday.}

\begin{abstract}
We prove that complete submanifolds, on which the Omori-Yau weak
maximum principle for the Hessian holds, with low codimension and bounded
by cylinders of small radius must have points rich in large positive extrinsic curvature. The lower the codimension is, the richer such points are. 
The smaller the radius is, the larger such curvatures are. This work 
unifies and generalizes several previous results on submanifolds with 
nonpositive extrinsic curvature.
\end{abstract}

\maketitle

\section{Introduction}\label{intrs}
The results of this article show that isometric immersions $f:M^m\to\tilde M^n$ with
low codimension and nonpositive extrinsic curvature at any point must satisfy strong
geometric conditions. The simplest result along this line is that a two-dimensional
surface with nonpositive curvature in $\R^3$ cannot be compact. This is a consequence of
the well-known fact that at a point of maximum of a distance function on a
compact surface in $\R^3$ the Gaussian curvature must be positive. It turns out that the
simple idea in the proof of this elementary fact has far-reaching generalizations for non
necessarily compact submanifolds in fairly general ambient Riemannian manifolds.

One of the main tools to extend this idea to higher dimensions and codimensions
is an algebraic lemma due to Otsuki \cite{MR0067550}, which will be stated in 
next section. On the other hand, a key ingredient to handle the noncompact
case is a maximum principle due to Omori \cite{MR0215259} and 
generalized by Pigola-Rigoli-Setti \cite{MR2116555}. 
Using this principle, Al\'ias-Bessa-Dajczer \cite{MR2529479} obtained an
estimate for the mean curvature of an isometric immersion 
$f:M^m\to N^{n+l}=P^n\times\R^l$, under some
assumptions on the manifold $P^n$, whose projection onto the first factor is bounded, the so-called cylindrically bounded submanifolds.
More recently, Al\'{i}as-Bessa-Montenegro \cite{MR2901222} have 
provided an estimate for the extrinsic curvatures  of such
submanifolds.

In the statement below and the sequel, $\rho$ stands for the distance function to a given reference point in $M^m$, $\log^{\left(j\right)}$ is the $j$-th iterate of the logarithm and $t\gg1$ means that $t$ is sufficiently large. Also $B_P[R]$ denotes the closed geodesic ball with radius $R$ centered at a point $o$ of $P^n$ and $\inj_P\left(o\right)$ is the injectivity radius of $P^n$ at $o$. Finally, $K_M\left(\sigma\right)$ denotes the sectional curvature of $M^m$ at a point $x\in M^m$ along the plane $\sigma\subset T_xM$, and similarly for $N^{n+l}$, $K_f\left(\sigma\right):=K_M\left(\sigma\right)-K_N\left(f_*\sigma\right)$ is the \emph{extrinsic sectional curvature} of $f$ at $x$ along $\sigma$ and $K^{\rad}_P$ is the radial sectional curvatures of $P^n$ with respect to $o$, that is, the sectional curvatures of tangent planes to $P^n$ containing the vector $\grad^Pr$, where $r$ is the distance function to $o$ in $P^n$.

\begin{theorem}[Al\'{i}as-Bessa-Montenegro, \cite{MR2901222}]\label{abmt}
Let $f:M^m\to N^{n+l}=P^n\times\R^l$ be an isometric immersion with codimension $p=n+l-m<m-l$ of a complete Riemannian manifold whose scalar curvature satisfies
\beeq\label{sca}
s_M\left(x\right)\geq -A^2\rho^2\left(x\right)\prod_{j=1}^{J}\left(\log^{\left(j\right)}\left(\rho\left(x\right)\right)\right)^2\virg\rho\left(x\right)\gg1,
\eneq
for some constant $A>0$ and some integer $J\geq1$. Assume that 
$f\left(M\right)\subset B_P[R]\times\R^l$, with $0<R<\min\left\{\inj_P\left(o\right),\frac{\pi}{2\sqrt{b}}\right\}$, where $\frac{\pi}{2\sqrt{b}}$ is replaced by $+\infty$ if $b\leq0$. If $K^{\rad}_P\leq b$ in $B_P[R]$, then
\beeq\label{abme}
\sup_M K_f\geq C_b^2\left(R\right),
\eneq
where
\begin{eqnarray*}
C_b(t)=\left\{\begin{array}{ll}
\sqrt{b}\cot(\sqrt{b}t)     & \If b>0\Aand 0<t<\dfrac{\pi}{2\sqrt{b}},\\
\frac{1}{t}                         & \If b=0\Aand t>0,\\
\sqrt{-b}\coth(\sqrt{-b}t)  & \If b<0\Aand t>0.
\end{array}\right.
\end{eqnarray*}
Moreover,
\beeq\label{abmie}
\sup_M K_M\geq C_b^2\left(R\right)+\inf_{B_P[R]}K_P.
\eneq
\end{theorem}

\bere
{\em
The geometry of the Euclidean factor $\R^l$ plays essentially no role in the proof of the above result. Indeed, estimate \eqref{abmie} remains true if the former is replaced by any Riemannian manifold $Q^l$, which need not be even complete, whereas for \eqref{abme} the only requirement is that $K_Q$ be bounded from above (see comment below Theorem \ref{svt}).}
\enre

It is worth pointing out that the codimension restriction $p<m-l$ cannot be relaxed. Actually, it implies that $n>2$ and $m>l+1$. In particular, in a three-dimensional ambient space $N^3$, that is, $n+l=3$, we have that $l=0$, and therefore 
$f\left(M\right)\subset B_P[R]$. In fact, the flat cylinder 
$\Sp^1\left(R\right)\times\R\subset B_{\R^2}[R]\times\R$ shows that the 
restriction $p<m-l$ is necessary.

On the other hand, estimates \eqref{abme} and \eqref{abmie} are sharp. Indeed, the function $C_b$ is well-known: the geodesic sphere $\partial B_{\Q_b^{m}}\left(R\right)$ of radius $R$ in the simply connected complete space form $\Q_b^{m}$ of constant sectional curvature $b$, with $R<\frac{\pi}{2\sqrt{b}}$ if $b>0$, is an umbilical hypersurface with principal curvatures being precisely $C_b(R)$. It shows that its extrinsic and intrinsic sectional curvatures are constant and equal to $C_b^2(R)$ and $C_b^2(R)+b$, respectively, the latter following from the former by the Gauss equation. Then, for every $n>2$ and $l\geq0$ we can consider 
$M^{m-1+l}=\partial B_{\Q_b^m}(R)\times\R^l$ and take 
$f:M^{m-1+l}\to B_{\Q_b^m}[R]\times\R^l$ to be the canonical isometric embedding. Therefore $\sup_MK_f$ and $\sup_MK_M$ are the constant extrinsic and intrinsic sectional curvatures $C_b^2(R)$ and $C_b^2(R)+b$ of 
$\partial B_{\Q_b^m}(R)$, respectively.

\vspace{.2cm}

The purpose of this paper is to provide a more accurate conclusion than the one of Theorem \ref{abmt} by precising how much extrinsic (respectively, intrinsic) sectional curvature satisfying estimate \eqref{abme} (respectively \eqref{abmie}) appears depending on how low the codimension is. The idea is that the lower the codimension is, the more extrinsic (respectively, intrinsic) sectional curvature satisfying \eqref{abme} (respectively \eqref{abmie}) will appear. 
Our main result can be stated as follows.

\begin{theorem}\label{svt}
Let $f:M^m\to N^{n+l}=P^n\times Q^l$ be an isometric immersion with 
codimension $p=n+l-m<m-l$ of a complete Riemannian manifold whose
radial sectional curvatures satisfy
\beeq\label{rsca}
K^{\rad}_M\left(x\right)\geq -A^2\rho^2\left(x\right)\prod_{j=1}^{J}\left(\log^{\left(j\right)}\left(\rho\left(x\right)\right)\right)^2\virg\rho\left(x\right)\gg1,
\eneq
for some constant $A>0$ and some integer $J\geq1$. Assume that 
$f(M)\subset B_P[R]\times Q^l$, with 
$0<R<\min\left\{\inj_P\left(o\right),\frac{\pi}{2\sqrt{b}}\right\}$, where $\frac{\pi}{2\sqrt{b}}$ is replaced by $+\infty$ if $b\leq0$. If $K^{\rad}_P\leq b$ in 
$B_P[R]$, then
\beeq\label{oe}
\sup_M\min\left\{\max_{\sigma\subset W}K_f\left(\sigma\right):\dim W>p+l\right\}\geq C_b^2\left(R\right).
\eneq
Moreover,
\beeq\label{oie}
\sup_M\min\left\{\max_{\sigma\subset W}K_M\left(\sigma\right):\dim W>p+l\right\}\geq C_b^2\left(R\right)+\inf_{B_P[R]}K_P.
\eneq
\end{theorem}

The estimates of Theorem \ref{svt} are clearly better than the ones of Theorem \ref{abmt}. Actually, \eqref{oe} and \eqref{oie} reduce to \eqref{abme} and \eqref{abmie}, respectively, only in the case of the highest allowed codimension $p=m-1-l$. On the other hand, although we make a stronger assumption on the curvature of $M^m$, if \eqref{sca} holds but \eqref{rsca} does not, then, since the scalar curvature is an average of sectional curvatures, we have that 
$\sup_MK_M=+\infty$, and hence \eqref{abmie} is trivially satisfied. Moreover, $K_P$ is clearly bounded in $B_P[R]$, thus if also $K_Q$ is bounded from above, we conclude that $\sup_MK_f=+\infty$ by the Gauss equation, so that \eqref{abme} also holds trivially in this case. Finally, note that the same example considered below Theorem \ref{abmt} shows that our estimates \eqref{oe} and \eqref{oie} are also sharp.

\bere\label{svtffmath}
{\em
Theorem \ref{svt} is a special case of the much stronger result, Theorem \ref{math}, given in Section \ref{proofs}.}
\enre

Let $f:M^m\to N^{n+l}=P^m\times Q^l$ be an isometric immersion. 
Recall that $f$ is said to be \emph{cylindrically bounded} if there exists
a (closed) geodesic ball $B_P[R]$ of $P^n$, centered at a point 
$o\in P^n$ with radius $R>0$, such that
\beeq\label{cbd}
f\left(M\right)\subset B_P[R]\times Q^l.
\eneq
Otherwise, we say that $f$ is \emph{cylindrically unbounded}. Denote by $R_f$ the \emph{extrinsic radius} of a cylindrically bounded isometric immersion $f$ (from $o$), that is, the smallest $R$ for which \eqref{cbd} holds. As a consequence of Theorem \ref{svt}, we have the following versions of the extrinsic radius results of Al\'{i}as-Bessa-Montenegro \cite{MR2901222}.

\becor\label{erthe}
Let $f:M^m\to N^{n+l}=P^n\times Q^l$ be an isometric immersion with 
codimension $p=n+l-m<m-l$ of a complete Riemannian manifold whose 
radial sectional curvatures satisfy \eqref{rsca}. Assume that $P^n$ is a 
complete Riemannian manifold with a pole and radial sectional curvatures 
$K^{\rad}_P\leq b\leq0$. If $f$ is cylindrically bounded, then 
\[
\sup_M\min\left\{\max_{\sigma\subset W}K_f\left(\sigma\right):
\dim W>p+l\right\}>-b
\]
and the extrinsic radius satisfies
\beeq
R_f\geq C_b^{-1}\left(\sqrt{\sup_M\min\left\{\max_{\sigma\subset W}K_f\left(\sigma\right):\dim W>p+l\right\}}\right).
\eneq
In particular, if 
\[
\sup_M\min\left\{\max_{\sigma\subset W}K_f\left(\sigma\right):
\dim W>p+l\right\}\leq-b,
\]
then $f$ is cylindrically unbounded.
\encor

\becor\label{erts}
Let $f:M^m\to N^{n+l}=\Sp^n\times Q^l$ be an isometric immersion with 
codimension $p=n+l-m<m-l$ of a complete Riemannian manifold whose
radial sectional curvatures satisfy \eqref{rsca}. If 
\[
\sup_M\min\left\{\max_{\sigma\subset W}K_M\left(\sigma\right):
\dim W>p+l\right\}\leq1,
\]
then
\beeq
R_f\geq\frac{\pi}{2}.
\eneq
\encor

On the other hand, a sharp lower bound for the Ricci curvature of bounded
complete Euclidean hypersurfaces was obtained by Leung \cite{MR707970}
and extended by Veeravalli \cite{MR1775899} to nonflat ambient space forms.

\begin{theorem}[Veeravalli, \cite{MR1775899}]\label{veer}
Let $f:M^n\to\Q_b^{n+1}$ be a complete hypersurface with sectional curvature bounded away from $-\infty$ such that 
$f\left(M\right)\subset B_{\Q_b^{n+1}}[R]$, with $R<\frac{\pi}{2\sqrt{b}}$ if $b>0$. Then
\beeq\label{estric}
\sup_M\Ric_M\geq C_b^2\left(R\right)+b,
\eneq
where $\Ric_M$ is the Ricci curvature of $M^n$.
\end{theorem}

Theorem \ref{svt} also gives an improvement of the above result, where we consider hypersurfaces of much more general ambient spaces and obtain that estimate \eqref{estric} actually holds for the scalar curvature. This shows the unifying character of our result.

\becor\label{clvsc}
Let $f:M^n\to P^{n+1}$ be a complete hypersurface whose radial sectional curvatures satisfy \eqref{rsca}. Assume that $f\left(M\right)\subset B_P[R]$,
with $R$ as in Theorem \ref{svt}. If $K^{\rad}_P\leq b$ in $B_P[R]$,
then
\beeqs
\sup_Ms_M\geq C_b^2\left(R\right)+\inf_{B_P[R]}K_P.
\eneqs
\encor

Again observe that for the geodesic sphere 
$M^n=\partial_{\Q_b^{n+1}}\left(R\right)$ of radius $R$ in $\Q_b^{n+1}$
the above inequality is in fact an equality. Corollary \ref{clvsc} leads to 
similar extrinsic radius results to Corollaries \ref{erthe} and \ref{erts} and,
in particular, a criterion of unboundness:

\becor
Let $f:M^n\to P^{n+1}$ be a complete hypersurface whose radial sectional curvatures satisfy \eqref{rsca}. Assume that $P^{n+1}$ is a complete 
Riemannian manifold with a pole and sectional curvatures $K_P\geq c$ and $K^{\rad}_P\leq b\leq0$. If $f(M)$ is bounded, then $\sup_Ms_M>c-b$ and
\beeq
R_f\geq C_b^{-1}\left(\sqrt{\sup_Ms_M-c}\right).
\eneq
In particular, if $\sup_Ms_M\leq c-b$, then $f\left(M\right)$ is unbounded.
\encor

\becor
Let $f:M^n\to\Sp^{n+1}$ be a complete hypersurface whose radial sectional curvatures satisfy \eqref{rsca}. If $\sup_Ms_M\leq1$, then
\beeq
R_f\geq\frac{\pi}{2}.
\eneq
\encor

Finally, we also generalize in the same spirit of Theorem \ref{svt} the 
second part of the work of Al\'{i}as-Bessa-Montenegro \cite{MR2901222}, concerning proper complete cylindrically bounded submanifolds with the 
norm of the second fundamental form with certain controlled growth.

\begin{theorem}\label{sffcth}
Let $f:M^m\to N^{n+l}=P^n\times Q^l$ be a proper isometric immersion
with codimension $p=n+l-m<m-l$ of a complete Riemannian manifold. 
Assume that $f\left(M\right)\subset B_P[R]\times Q^l$, with $R$ as in 
Theorem \ref{svt}. If $K^{\rad}_P\leq b$ in $B_P[R]$, $Q^l$ is a complete Riemannian manifold with a pole and
\beeq\label{bsff}
\sup_{f^{-1}(B_P[R]\times\partial B_Q\left(t\right))}\left\|\al\right\|
\leq\varsigma\left(t\right),
\eneq
where $\al$ is the second fundamental form of $f$ and 
$\varsigma:\left[0,+\infty\right)\to\R$ is a positive function satisfying 
$\int_0^{+\infty}1/\varsigma=+\infty$, then \eqref{oe} and \eqref{oie} hold.
\end{theorem}

For hypersurfaces, the growth rate of the norm of the second fundamental
form can be improved as follows.

\begin{theorem}\label{ltfh}
Let $f:M^n\to N^{n+1}=P^{n+1-l}\times Q^l\virg n-l>1$, be a properly immersed complete hypersurface with $f\left(M\right)\subset B_P[R]\times Q^l$, 
with $R$ as in Theorem \ref{svt}. Suppose that $N^{n+1}$ satisfies the assumptions as in Theorem \ref{sffcth} and the second fundamental form
$\al$ satisfies
\beeq\label{hyphyp}
\sup_{f^{-1}\left(B_P[R]\times\partial B_Q\left(t\right)\right)}\left\|\al\right\|\leq\varsigma^2\left(t\right),
\eneq
where $\varsigma:\left[0,+\infty\right)\to\R$ is a positive function satisfying
\beeqs
\int_0^{+\infty}\frac{1}{\varsigma}=+\infty\Aand\limsup_{t\to+\infty}\frac{1}{\varsigma\left(t\right)}<+\infty.
\eneqs
Then \eqref{oe} and \eqref{oie} hold.
\end{theorem}

\section{Preliminaries}

Our main tools to build the proof of Theorem \ref{svt} are Otsuki's Lemma, 
the Omori-Yau maximum principle and the Hessian comparison theorem, 
which for the sake of organization will be presented in two subsections.

\subsection{Otsuki's Lemma}

Throughout this subsection, $V^n$ and $W^p$ will be real vector spaces
of dimensions $n$ and $p$, respectively, endowed with positive definite 
inner products. For a symmetric bilinear form $\al:V^n\times V^n\to W^p$,
we denote
\beeqs
K_{\al}\left(X,Y\right)=\langle\al\left(X,X\right),\al\left(Y,Y\right)\rangle-
\left\|\al\left(X,Y\right)\right\|^2,
\eneqs
for any pair of vectors $X,\,Y\in V^n$. If $\sigma$ is a two-dimensional 
subspace of $V^n$, we define
\beeqs
K_{\al}\left(\sigma\right)=\frac{K_{\al}\left(X,Y\right)}{\left\|X\wedge Y\right\|^2},
\eneqs
where $\left\{X,Y\right\}$ is any basis of $\sigma$ and 
$\left\|X\wedge Y\right\|^2=\left\|X\right\|^2\left\|Y\right\|^2-\langle X,Y\rangle^2$.
Given an isometric immersion $f:M^m\to\tilde M^n$ with second fundamental
form $\al$, then for any $x\in M^m$ and any plane $\sigma\subset T_xM$
the Gauss equation yields
\beeqs
K_{\al}\left(\sigma\right)=K_f\left(\sigma\right).
\eneqs
A basic tool in this article is the following algebraic lemma, known as Otsuki's Lemma (for a
proof see, for instance, \cite{MR1393941}).
\bele\label{olwv}
Let $\al:V^n\times V^n\to W^p$ be a symmetric bilinear form. Suppose there exists a real number $\lambda\geq0$ such that
\begin{enumerate}[(i)]
\item $K_{\al}\left(\sigma\right)\leq\lambda$ for every plane $\sigma\subset V^n$,
\item $\left\|\al\left(X,X\right)\right\|>\sqrt{\lambda}$ for every unit vector $X\in V^n$.
\end{enumerate}
Then $p\geq n$.
\enle

Given a symmetric bilinear form $\al:V^n\times V^n\to W^p$, a vector 
$X\in V^n$ is said to be an \emph{asymptotic vector} of $\al$ if
\beeqs
\al\left(X,X\right)=0.
\eneqs
In the next statement and the sequel, we write $K_{\al}\leq0$ (respectively, $K_{\al}<0$) as a shorthand for $K_{\al}\left(\sigma\right)\leq0$ (respectively, $K_{\al}\left(\sigma\right)<0$) for any plane $\sigma\subset V^n$.

\becor\label{colwv}
Let $\al:V^n\times V^n\to W^p$ be a symmetric bilinear form.
\begin{enumerate}[(i)]
\item If $K_{\al}\leq0$, then any subspace $S\subset V^n$, with 
$\dim S>p$, contains a nonzero asymptotic vector of $\al$.
\item If $K_{\al}<0$, then $p\geq n-1$.
\end{enumerate}
\encor
\bepr
(\emph{i}) This is just an equivalent way of stating Lemma \ref{olwv} for $\lambda=0$.\vspace{4pt}

\noindent (\emph{ii}) If there are no nonzero asymptotic vectors of $\al$, the result follows from Lemma \ref{olwv}. Suppose $p<n-1$, and assume that there exists a nonzero vector $X_0\in V^n$ such that $\al\left(X_0,X_0\right)=0$. Denote by $U$ the orthogonal complement to $X_0$ in $V^n$, and consider the
linear map $B_{X_0}:U\to W^p$ defined by $B_{X_0}\left(Y\right)=\al\left(X_0,Y\right)$. Since $\dim U=n-1>p$,
there exists a nonzero vector $Y_0\in U$ such that $B_{X_0}\left(Y_0\right)=0$. This fact, together with
$\al\left(X_0,X_0\right)=0$, contradicts the assumption.
\enpr

The following result is a direct consequence of Corollary \ref{colwv}-(\emph{ii}).

\begin{theorem}\label{tols}
Let $f:M^n\to\tilde M^{n+p}$ be an isometric immersion. Assume that there
exists a point $x_0\in M^n$ and a subspace $V_{x_0}\subset T_{x_0}M$ 
with dimension $d$ such that $K_f\left(\sigma\right)<0$ along every plane 
$\sigma\subset V_{x_0}$. Then $p\geq d-1$.
\end{theorem}

%\bere
The preceding inequality is sharp, as shown, e.g., by the $n$-dimensional Clifford
torus in $\Sp^{2n-1}$. %\enre
Theorem \ref{tols} comes from a purely algebraic restriction on the ``codimension'' $p$ of symmetric bilinear forms $\al:V^n\times V^n\to W^p$ with $K_{\al}<0$, which gives its punctual nature. If, on the other hand, $K_f\left(\sigma\right)\leq0$ in the above statement, it is possible to use part (\emph{i}) of Corollary \ref{colwv} to obtain the stronger restriction $p\geq d$, provided that the subspace $V_{x_0}$ is free of asymptotic directions. Actually, this is a central idea in the study of the global implications of nonpositive extrinsic curvature in low codimension. The presence of a certain amount of nonpositive extrinsic sectional curvature everywhere together with some global assumption that can guarantee the existence of points without asymptotic directions must imply codimension restrictions. For compact Riemannian manifolds, for instance, one obtains the following generalization of a result for the flat case due to Tompkins \cite{MR1546106}.
\begin{theorem}
Let $M^n$ be a compact Riemannian manifold such that at any point $x\in M^n$ there exists a subspace $V_x$ of $T_xM$ with dimension $d$ such that $K_M\left(\sigma\right)\leq0$ for every plane $\sigma\subset V_x$. If $f:M^n\to\R^{n+p}$ is an isometric immersion, then $p\geq d$.
\end{theorem}
\bepr
Since $M^n$ is compact, it is well known that there exist a point $x_0\in M^n$ and a normal vector $\xi\in N_fM\left(x_0\right)$ such that the shape operator $A_\xi$ is positive definite, and in particular $\al_{x_0}\left(X,X\right)\neq0$ for every nonzero vector $X\in T_{x_0}M$. Furthermore, $K_{\al_{x_0}}\left(\sigma\right)\leq0$ for every
plane $\sigma\subset V_{x_0}$ by the Gauss equation. The statement then follows from Corollary \ref{colwv}-(\emph{i}).
\enpr

For the noncompact case, on the other hand, we have the following immediate consequence of our Corollary \ref{erthe}.

\becor\label{corotstom}
Let $f:M^n\to P^{n+p}$ be an isometric immersion of a
complete Riemannian manifold whose radial sectional curvatures satisfy \eqref{rsca} into a Hadamard manifold. Assume that at any point $x\in M^n$ there exists a subspace $V_x$ of $T_xM$ with dimension $d$ such that $K_f\left(\sigma\right)\leq0$ for every plane $\sigma\subset V_x$. If $p<d$, then $f\left(M\right)$ is unbounded.
\encor

\subsection{Omori-Yau maximum principle and Hessian comparison theorem}

The \emph{Omori-Yau maximum principle for the Hessian} is said to hold on a given
Riemannian manifold $M^n$ if, for any function $g\in C^2\left(M\right)$ with $g^*=\sup_Mg<+\infty$,
there exists a sequence of points $\left\{x_k\right\}_{k\in\N}$ in $M^n$ satisfying:
\begin{enumerate}
\item[(i)] $g\left(x_k\right)>g^*-\frac{1}{k}$,
\item[(ii)] $\left\|\grad^Mg\left(x_k\right)\right\|<\frac{1}{k}$,
\item[(iii)] $\hess^Mg\left(x_k\right)\left(X,X\right)\leq\frac{1}{k}\langle X,X\rangle$ for all $X\in T_{x_k}M$.
\end{enumerate}
Such a sequence $\left\{x_k\right\}_{k\in\N}$ satisfying (\emph{i})-(\emph{iii})
above is called an \emph{Omori-Yau Hessian sequence} for $g$. 
One says that the \emph{Omori-Yau maximum principle} holds if the 
condition on the
Hessian is replaced by a similar one on the Laplacian, namely, if
\begin{enumerate}
\item[(iii)] $\Delta^Mg\left(x_k\right)\leq\frac{1}{k}$.
\end{enumerate}
In this case, $\left\{x_k\right\}_{k\in\N}$ is called an 
\emph{Omori-Yau sequence} for $g$.

\vspace{.2cm}

The following is a function theoretic characterization of Riemannian manifolds
that satisfy the Omori-Yau maximum principle for the Hessian. For the proof of this, as well as of the other results related to the Omori-Yau maximum principle in this subsection, we refer to \cite{MR2116555}.

\begin{theorem}\label{oympct}
Assume that the Riemannian manifold $M^n$ supports a nonnegative function $\gamma\in C^2\left(M\right)$ satisfying the following conditions:
\begin{enumerate}[(a)]
\item $\gamma$ is proper, that is, $\gamma\left(x\right)\to+\infty$ as $x\to\infty$,
\item $\left\|\grad^M\gamma\right\|\leq c\sqrt{\gamma}$ outside a compact subset of $M^n$ for some constant $c>0$,
\item $\hess^M\gamma\leq d\sqrt{\gamma F\left(\sqrt{\gamma}\right)}\langle\,,\,\rangle$ outside a compact subset of $M^n$, for some $d>0$ and some $F\in C^\infty\left(\left[0,+\infty\right)\right)$ that satisfies:
\beeqs
\mbox{(i) }F\left(0\right)>0,\mbox{ (ii) }F'\left(t\right)\geq0\mbox{ on }\left[0,+\infty\right),\mbox{ (iii) }1/\sqrt{F\left(t\right)}\notin L^1\left[0,+\infty\right).
\eneqs
\end{enumerate}
Then, the Omori-Yau maximum principle for the Hessian holds on $M^n$.
\end{theorem}

We point out that, although completeness of $M^n$ is not required in 
Theorem \ref{oympct}, it is a consequence of the assumptions (\emph{a})
and (\emph{b}). Examples of functions satisfying
the requirements in Theorem \ref{oympct} are given by
\beeqs
F\left(t\right)=A^2t^2\prod_{j=1}^J\left(\log^{\left(j\right)}t\right)^2,\,t\gg1,
\eneqs
where $A>0$ is a constant and $J\geq1$ is an integer.

\vspace{.2cm}

Sometimes, for the applications of the maximum principle as in our paper,
the following weaker version is enough.

\vspace{.2cm}

The \emph{Omori-Yau weak maximum principle for the Hessian} is said to hold on a
Riemannian manifold $M^n$ if for any function $g\in C^2\left(M\right)$ with $g^*=\sup_Mg<+\infty$
there exists a sequence of points $\left\{x_k\right\}_{k\in\N}$ satisfying:
\begin{enumerate}[(i)]
\item $g\left(x_k\right)>g^*-\frac{1}{k}$,
\item $\hess^Mg\left(x_k\right)\left(X,X\right)\leq\frac{1}{k}\langle X,X\rangle$ for all $X\in T_{x_k}M$.
\end{enumerate}
A sequence $\left\{x_k\right\}_{k\in\N}$ satisfying (\emph{i}) and (\emph{ii}) is called an \emph{Omori-Yau weak Hessian sequence} for $g$.

\vspace{.2cm}

Riemannian manifolds that satisfy the Omori-Yau weak maximum principle for
the Hessian are characterized as follows.

\begin{theorem}\label{oywmpct}
Assume that the Riemannian manifold $M^n$ supports a nonnegative function $\gamma\in C^2\left(M\right)$ satisfying the following conditions:
\begin{enumerate}[(a)]
\item $\gamma$ is proper,
\item $\hess^M\gamma\leq d\sqrt{\gamma F\left(\sqrt{\gamma}\right)}\langle\,,\,\rangle$ outside a compact subset of $M^n$, for some $d>0$ and some $F\in C^\infty\left(\left[0,+\infty\right)\right)$ as in Theorem \ref{oympct}.
\end{enumerate}
Then, the Omori-Yau weak maximum principle for the Hessian holds on $M^n$.
\end{theorem}

Accordingly, it is said that the \emph{Omori-Yau weak maximum principle} holds if (\emph{ii})
is replaced by the condition
\begin{itemize}
\item[(ii)] $\Delta^Mg\left(x_k\right)\leq\frac{1}{k}$,
\end{itemize}
in which case $\left\{x_k\right\}_{k\in\N}$ is called an \emph{Omori-Yau weak sequence} for $g$.

\vspace{.2cm}

The function theoretic approach to the Omori-Yau Maximum Principle given
in Theorem \ref{oympct} allows one to apply it in different situations, where
the choices of $\gamma$ and $F$ are suggested by the geometric setting. 
For instance, one has the following special case, where as previously 
agreed $\rho$ stands for the distance function on a Riemannian manifold 
$M^n$ to a fixed point.

\begin{theorem}\label{oympsc}
Let $M^n$ be a complete noncompact Riemannian manifold. Assume that 
$K^{\rad}_M\geq-F\left(\rho\right)$, where $F\in C^\infty\left[0,+\infty\right)$ 
satisfies the conditions listed in Theorem \ref{oympct} and is even at the 
origin, that is, its derivatives satisfy $F^{\left(2k+1\right)}\left(0\right)=0$ 
for $k\geq0$. Then, the Omori-Yau maximum principle for the Hessian 
holds on $M^n$.
\end{theorem}

\bere
{\em
If we only assume that $\Ric_M\left(\grad^M\rho\right)\geq-F\left(\rho\right)$,
then the conclusion is that the Omori-Yau maximum principle holds on 
$M^n$.}
\enre

The last ingredient for the proof of Theorem \ref{svt} is the following version
of the well-known Hessian comparison theorem given in \cite{MR1333601}.

\begin{theorem}\label{hcth}
Let $M^n$ be a Riemannian manifold and $o,x\in M^n$ be such that there
is a minimizing geodesic $\gamma$ joining $o$ and $x$, and let $\rho$ be
the distance function to $o$. Suppose that $K^{\rad}_M\leq b$ along 
$\gamma$. If $b>0$ assume $\rho\left(x\right)<\frac{\pi}{2\sqrt{b}}$. Then,
we have
\beeqs
\hess^M\rho\geq C_b\left(\rho\right)\left(\langle\,,\,\rangle-d\rho\otimes d\rho\right)
\eneqs
along $\gamma$.
\end{theorem}

\section{Proofs}\label{proofs}
Let $f:M^m\to\tilde M^n$ be an isometric immersion between Riemannian 
manifolds. Given a function $h\in C^\infty(\tilde M)$ we set 
$g=h\circ f\in C^\infty(M)$. Since
\beeqs
\langle\grad^Mg(x),X\rangle=\langle\grad^{\tilde M}h(f(x)),f_*X\rangle
\eneqs
for every $x\in M^n\Aand X\in T_xM$, we obtain
\beeq\label{gradcomp}
f_*\grad^Mg(x)=\left(\grad^{\tilde M}h(f(x))\right)^T,
\eneq
where $(\,)^T$ is the tangent component. An easy computation using
the Gauss formula gives the well-known relation 
(see e.g. \cite{MR623135}):
\beeq\label{hc}
\hess^Mg(x)(X,Y)=\hess^{\tilde M}h(f(x))(f_*X,f_*Y)+
\langle\grad^{\tilde M}h(f(x)),\al_x(X,Y)\rangle
\eneq
for all $x\in M^n\Aand X,\,Y\in T_xM$, where $\al_x$ stands for the second fundamental form of $f$ at $x$. In particular, taking traces with respect to 
an orthonormal frame $\left\{e_1,\dots,e_m\right\}$
in $T_xM$ yields
\beeqs
\Delta^Mg(x)=\sum_{i=1}^m\hess^{\tilde M}h(f(x))(f_*e_i,f_*e_i)
+n\langle\grad^{\tilde M}h(f(x)),H(x)\rangle,
\eneqs
where $H(x)=\frac{1}{n}\sum_{i=1}^m\al_x(e_i,e_i)$ is the mean 
curvature vector of $f$ at $x$.

\vspace{.2cm}

Given an isometric immersion $f:M^m\to N^{n+l}=P^n\times Q^l$, we
denote by $\pi_P:N^{n+l}\to P^n$ and $\pi_Q:N^{n+l}\to Q^l$ the 
projections onto $P^n$ and $Q^l$, respectively. We write $(y,z)$
for points in $N^{n+l}=P^n\times Q^l$ and by abuse of notation 
denote $y=\pi_P\circ f\Aand z=\pi_Q\circ f$.

Moreover, set
\begin{eqnarray*}
\psi_b=\left\{\begin{array}{ll}
1-\cos\left(\sqrt{b}t\right)     & \If b>0,\\
t^2                        & \If b=0,\\
\cosh\left(\sqrt{-b}t\right)  & \If b<0,
\end{array}\right.
\end{eqnarray*}
where $t>0\If b\leq0\Aand 0<t<\frac{\pi}{2\sqrt{b}}\If b>0$. Hence
$\psi_b''=C_b\psi_b'$. We define $h\in C^\infty\left(N\right)$ by 
$h=\psi_b\circ r\circ\pi_P$, where $r$ is the distance function on $P^n$ 
to the reference point $o$. We call $g=h\circ f$ the 
\emph{modified radial function} of $f$.

\subsection{Proofs of Theorem \ref{svt} and corollaries}
As mentioned in Remark \ref{svtffmath}, Theorem \ref{svt} is a consequence
of the following stronger result.  Here the \emph{algebraic codimension}
$p(x)$ of an isometric immersion $f:M^m\to\tilde M^n$ at $x\in M^m$
is the dimension of its first normal space $N_1(x)$ and a sequence 
of real numbers $\left\{p_k\right\}_{k\in\N}$ is said to be \emph{strictly 
bounded from above} by another $\left\{q_k\right\}_{k\in\N}$ if 
$p_k<q_k$ for all $k\in\N$.

\begin{theorem}\label{math}
Let $f:M^m\to N^{n+l}=P^n\times Q^l$ be an isometric immersion. Assume
that $f\left(M\right)\subset B_P[R]\times Q^l$, with $0<R<\min\left\{\inj_P\left(o\right),\frac{\pi}{2\sqrt{b}}\right\}$, where $\frac{\pi}{2\sqrt{b}}$ is replaced by $+\infty$ if $b\leq0$. If $K^{\rad}_P\leq b$ in $B_P[R]$, then 
\beeq\label{lie}
\liminf\min\left\{\max_{\sigma\subset W}K_f\left(\sigma\right):
\dim f_*W\cap T_{y\left(x_k\right)}P>p(x_k)\right\}
\geq C_b^2\left(R\right)
\eneq
for all Omori-Yau weak Hessian sequence $\left\{x_k\right\}_{k\in\mathbb{N}}$ 
for the modified radial function of $f$ with algebraic codimension sequence
$\left\{p\left(x_k\right)\right\}_{k\in\mathbb{N}}$ strictly bounded from above
by $\left\{\dim f_*T_{x_k}M\cap T_{y\left(x_k\right)}P\right\}_{k\in\mathbb{N}}$. Furthermore,
\beeq\label{liie}
\liminf\min\left\{\max_{\sigma\subset W}K_M\left(\sigma\right):
\dim f_*W\cap T_{y\left(x_k\right)}P>p(x_k)\right\}
\geq C_b^2\left(R\right)+\inf_{B_P[R]}K_P.
\eneq
\end{theorem}
\bepr
By the assumption that $f\left(M\right)\subset B_P[R]\times Q^l$, 
the modified radial function $g$ satisfies
\beeqs
g^*\leq\psi_b\left(R\right)<+\infty,
\eneqs
where we write $\left(\,\right)^*=\sup_M\left(\,\right)$. Let $\left\{x_k\right\}_{k\in\mathbb{N}}$ be an Omori-Yau weak Hessian sequence for $g$, that is,
\begin{itemize}
\item[(i)] $g\left(x_k\right)>g^*-\frac{1}{k}$,
\item[(ii)] $\hess^Mg\left(x_k\right)\leq\frac{1}{k}\langle\hspace{1pt},\hspace{1pt}\rangle.$
\end{itemize}
For each $k\in\N$, take a tangent subspace $W_k\subset T_{x_k}M$ such that $\dim V_k>p_k$, where $V_k=f_*^{-1}\left(f_*W_k\cap T_{y_k}P\right)$ and for simplicity of notation we write $p_k=p\left(x_k\right)$, $y_k=y\left(x_k\right)$. If the algebraic codimension sequence $\left\{p_k\right\}_{k\in\N}$ is strictly bounded from above by $\left\{\dim f_*T_{x_k}M\cap T_{y_k}P\right\}_{k\in\mathbb{N}}$, then at least $W_k=T_{x_k}M$ satisfies the latter condition, so that the sets on the left-hand side of inequalities \eqref{lie} and \eqref{liie} are nonempty. The idea of the argument is to use (\emph{ii}) above and \eqref{hc} to estimate $\left\|\al_{x_k}\left(X,X\right)\right\|$ for $X\in V_k$, and then apply Lemma \ref{olwv} to 
$\al_{x_k}|_{V_k\times V_k}$. This will imply the estimates in the statement.
By \eqref{gradcomp}, we have
\beeqs
\grad^Nh\left(f\left(x\right)\right)=f_*\grad^Mg\left(x\right)+\left(\grad^Nh\left(f\left(x\right)\right)\right)^\perp,
\eneqs
where $\left(\,\right)^\perp$ denotes the normal component. Note that
\beeq\label{gmrf}
\grad^Nh\left(y,z\right)=\psi_b'\left(r\left(y\right)\right)\grad^Pr\left(y\right).
\eneq
Since $h$ only depends on $P^n$, we obtain from \eqref{hc} and \eqref{gmrf} that
\begin{multline*}
\hess^Mg\left(x\right)\left(X,X\right)=\hess^P\psi_b\circ r\left(y\left(x\right)\right)\left(X_P,X_P\right)\\
+\psi_b'\left(r\left(y\left(x\right)\right)\right)\langle\grad^Pr\left(y\left(x\right)\right),\al_x\left(X,X\right)\rangle,
\end{multline*}
where $X_P=y_*X$. Observe that
%\footnotesize
%\beeqs
\begin{multline*}
\hess^P\psi_b\circ r\left(y\right)\left(X_P,X_P\right)=
\psi_b''\left(r\left(y\right)\right)\langle\grad^Pr\left(y\right),X_P\rangle^2 \\ 
+\psi_b'\left(r\left(y\right)\right)\hess^Pr\left(y\right)\left(X_P,X_P\right).
%\eneqs
\end{multline*}
%\normalsize
Since $\psi_b''=C_b\psi_b'$, the last two equations yield
\begin{multline}\label{mrfp}
\hess^Mg\left(x\right)\left(X,X\right)=\psi_b'\left(r\left(y\left(x\right)\right)\right)(C_b\left(r\left(y\left(x\right)\right)\right)\langle\grad^Pr\left(y\left(x\right)\right),X_P\rangle^2\\
+\langle\grad^Pr\left(y\left(x\right)\right),\al_x\left(X,X\right)\rangle+\hess^Pr\left(y\left(x\right)\right)\left(X_P,X_P\right)).
\end{multline}
Theorem \ref{hcth} gives
\begin{eqnarray}\label{hct}
\begin{aligned}
\hess^Pr\left(y\right)\left(Y,Y\right) &=
\hess^Pr\left(y\right)\left(Y^\perp,Y^\perp\right) \\
&\geq C_b\left(r\left(y\right)\right)\left(\left\|Y\right\|^2-\langle\grad^Pr\left(y\right),Y\rangle^2\right),
\end{aligned}
\end{eqnarray}
where $Y\in T_yP$ and here $Y^\perp$ is defined by the orthogonal decomposition
\beeqs
Y=\langle\grad^Pr\left(y\right),Y\rangle\grad^Pr\left(y\right)+Y^\perp.
\eneqs
Now, since $X_P=f_*X$ for any $X\in V_k$, we obtain from \eqref{mrfp} and \eqref{hct} that
\beeqs
\begin{aligned}
\hess^Mg\left(x_k\right)\left(X,X\right) & \geq\psi_b'\left(r_k\right)\left(C_b\left(r_k\right)\left\|X\right\|^2+\langle\grad^Pr\left(y_k\right),\al_{x_k}\left(X,X\right)\rangle\right)\\
                   & \geq\psi_b'\left(r_k\right)\left(C_b\left(r_k\right)\left\|X\right\|^2-\left\|\al_{x_k}\left(X,X\right)\right\|\right),
\end{aligned}
\eneqs
where $r_k=r\left(y_k\right)$. Hence, by (\emph{ii})
\beeqs
\frac{1}{k}\left\|X\right\|^2\geq\psi_b'\left(r_k\right)\left(C_b\left(r_k\right)\left\|X\right\|^2-\left\|\al_{x_k}\left(X,X\right)\right\|\right)
\eneqs
for every $x_k$ and every $X\in V_k$, and therefore,
\beeqs
\left\|\al_{x_k}\left(X,X\right)\right\|\geq\left(C_b\left(r_k\right)-\frac{1}{k\psi_b'\left(r_k\right)}\right)\left\|X\right\|^2.
\eneqs
Since $g\left(x_k\right)=\psi_b\left(r_k\right)$ approaches $g^*=\psi_b\left(r^*\right)>0$ by (\emph{i}) and $\psi_b|_{\left[0,R\right]}$ is a homeomorphism onto its image (recall that $R<\frac{\pi}{2\sqrt{b}}$ if $b>0$), it follows that $r_k$ goes to $r^*>0$, and in particular $\psi_b'\left(r_k\right)\to\psi_b'\left(r^*\right)>0$. Thus,
\beeqs
C_b\left(r_k\right)-\frac{1}{k\psi_b'\left(r_k\right)}>0
\eneqs
for $k$ sufficiently large and, as $\dim V_k>p_k$, we can apply Lemma \ref{olwv} to
\beeqs
\al_{x_k}|_{V_k\times V_k}:V_k\times V_k\to N_1\left(x_k\right).
\eneqs
We obtain a plane $\sigma_k\subset V_k$ such that, by the Gauss equation,
\beeqs
K_f\left(\sigma_k\right)=K_{\al_{x_k}}\left(\sigma_k\right)\geq\left(C_b\left(r_k\right)-\frac{1}{k\psi_b'\left(r_k\right)}\right)^2.
\eneqs
In particular,
\beeqs
\max_{\sigma\subset W_k}K_f\left(\sigma\right)\geq\left(C_b\left(r_k\right)-\frac{1}{k\psi_b'\left(r_k\right)}\right)^2,
\eneqs
but since the subspaces $W_k\subset T_{x_k}M$ satisfying $\dim f_*\left(W_k\right)\cap T_{y_k}P>p_k$ have been taken arbitrarily, we have indeed
\beeqs
\min\left\{\max_{\sigma\subset W_k}K_f\left(\sigma\right):\dim f_*W_k\cap T_{y_k}P>p_k\right\}\geq\left(C_b\left(r_k\right)-\frac{1}{k\psi_b'\left(r_k\right)}\right)^2.
\eneqs
Then \eqref{lie} follows by letting $k\to+\infty$. We will now compare the sectional curvatures $K_M\left(\sigma_k\right)\Aand K_N\left(f_*\sigma_k\right)$. Since $\sigma_k\subset V_k\perp T_{z_k}Q,\ z_k=z\left(x_k\right)$, then
\beeqs
K_N\left(f_*\sigma_k\right)=K_P\left(y_*\sigma_k\right).
\eneqs
Then, we have that
\beeqs
\begin{aligned}
K_M\left(\sigma_k\right)&=K_f\left(\sigma_k\right)+K_P\left(y_*\sigma_k\right)\\
&\geq\left(C_b\left(r_k\right)-\frac{1}{k\psi_b'\left(r_k\right)}\right)^2+\inf_{B_P\left[R\right]}K_P,
\end{aligned}
\eneqs
and \eqref{liie} follows by a similar argument.
\enpr

\bere
{\em
That the maximum and minimum on the left-hand side of \eqref{lie} (and similarly for \eqref{liie}) are indeed attained can be argued as follows. At each $x=x_k$, the extrinsic sectional curvature $K_f=K_{\al_x}:G_2\left(T_xM\right)\to\R$ is a continuous function on the Grassmannian $G_2\left(T_xM\right)$ of (nonoriented) planes in $T_xM$, and by compactness attains its maximum and minimum. Since $G_2\left(W\right)$ is a compact subset of $G_2\left(T_xM\right)$ for any subspace $W\subset T_xM$, so does the restriction $K_f|_W$. Let 
\[
\left\{W_j\right\}_{j\in\N}\subset\mathcal{W}:=
\left\{W\subset T_xM:\dim f_*W\cap T_{y\left(x\right)}P>p
\left(x\right)\right\}\neq\emptyset
\]
be a sequence such that
\beeqs
\max_{\sigma\subset W_j}K_f\left(\sigma\right)\to\inf\left\{\max_{\sigma\subset W}K_f\left(\sigma\right):W\in\mathcal{W}\right\},
\eneqs
as $j\to+\infty$. After passing to a subsequence we can without loss of generality assume that all $W_j$ have the same dimension $d$ and converge to some $W_\infty\in G_d\left(T_xM\right)$, where $G_d\left(T_xM\right)$ denotes the Grassmannian of (nonoriented) $d$-planes in $T_xM$. Moreover, since $f_*$ is an isomorphism onto its image, it is clear that the function $W\in G_d\left(T_xM\right)\mapsto \dim f_*W\cap T_{y\left(x\right)}P$ is upper semicontinuous, and so
\beeqs
\dim f_*\left(W_\infty\right)\cap T_{y\left(x\right)}P>p\left(x\right),
\eneqs
or equivalently, $W_\infty\in\mathcal{W}$. Hence, 
\[
\max_{\sigma\subset W_\infty}K_f\left(\sigma\right)\in
\left\{\max_{\sigma\subset W}K_f\left(\sigma\right):W\in\mathcal{W}\right\}.
\]
Finally, a straightforward contradiction argument allows to conclude that
\beeqs
\lim_{j\to+\infty}\max_{\sigma\subset W_j}K_f\left(\sigma\right)=\max_{\sigma\subset W_\infty}K_f\left(\sigma\right),
\eneqs
and therefore 
\[
\max_{\sigma\subset W_\infty}K_f\left(\sigma\right)=
\min\left\{\max_{\sigma\subset W}K_f\left(\sigma\right):W\in\mathcal{W}\right\}.
\]}
\enre

\bepr[Proof of Theorem \ref{svt}]
According to Theorem \ref{oympsc}, the curvature decay in the statement is sufficient to
conclude that the Omori-Yau maximum principle for the Hessian holds 
on $M^m$. Thus there exists in $M^m$ an Omori-Yau Hessian sequence
$\left\{x_k\right\}_{k\in\N}$ for the modified radial function of $f$, whose
algebraic codimension sequence $\left\{p_k\right\}_{k\in\N}$ satisfies
\beeqs
p_k\leq p<m-l\leq\dim f_*T_{x_k}M\cap T_{y_k}P,
\eneqs
that is, $\left\{p_k\right\}_{k\in\N}$ is strictly bounded from above by $\left\{\dim f_*T_{x_k}M\cap T_{y_k}P\right\}$ as required in Theorem \ref{math}. Moreover, given a subspace $W\subset T_{x_k}M$ with $\dim W>p+l$, it is clear that
$\dim f_*W\cap T_{y_k}P>p\geq p_k$. In other words,
\beeqs
\left\{W\subset T_{x_k}M:\dim W>p+l\right\}\subset\left\{W\subset T_{x_k}M:\dim f_*W\cap T_{y_k}P>p_k\right\},
\eneqs
and in particular
\beeqs
\left\{\max_{\sigma\subset W}K_f\left(\sigma\right):\dim W>p+l\right\}\subset\left\{\max_{\sigma\subset W}K_f\left(\sigma\right):\dim f_*W\cap T_{y_k}P>p_k\right\}.
\eneqs
Therefore,
%\small
\beeqs
\min\left\{\max_{\sigma\subset W}K_f\left(\sigma\right):\dim W>p+l\right\}\geq\min\left\{\max_{\sigma\subset W}K_f\left(\sigma\right):\dim f_*W\cap T_{y_k}P>p_k\right\},
\eneqs
\normalsize
and \eqref{oe} follows immediately from \eqref{lie}. Similarly, \eqref{oie} 
follows from \eqref{liie}. 
\enpr

\bepr[Proof of Corollary \ref{erthe}]
Follows immediately from \eqref{oe} observing that $\inf C_b=\sqrt{-b}$ for $b\leq0$.
\enpr

\bere
{\em
If, in addition, $K_P\geq c$ in Corollary \ref{erthe}, then \eqref{oie} provides
\beeqs
\sup_M\min\left\{\max_{\sigma\subset W}K_M\left(\sigma\right):\dim W>p+l\right\}\leq c-b
\eneqs
as a criterion of cylindrical unboundness. In particular, when $P^{n+p}$
is either the Euclidean space $\R^{n+p}$ or the hyperbolic space 
$\mathbb{H}^{n+p}$ and $M^n$ is a complete Riemannian manifold 
(whose radial sectional curvatures satisfy \eqref{rsca}) in which at any 
point $x\in M^n$ there exists a subspace $V_x$ of $T_xM$ with dimension
$d$ such that $K_M(\sigma)\leq0$ for every plane 
$\sigma\subset V_x$, we conclude that every isometric immersion 
$f:M^n\to P^{n+p}$ with codimension $p<d$ is unbounded (compare 
with Corollary \ref{corotstom}).}
\enre

\bepr[Proof of Corollary \ref{clvsc}]
Here $p=1$ and $l=0$, so that \eqref{oie} in Theorem \ref{svt} yields
\beeqs
\sup_M\min K_M\geq C_b^2\left(R\right)+\inf_{B_P[R]}K_P.
\eneqs
Since clearly $s_M\geq\min K_M$ at each point, the corollary follows.
\enpr

\subsection{Proofs of Theorems \ref{sffcth} and \ref{ltfh}}
\bepr[Proof of Theorem \ref{sffcth}] Again consider the modified radial function $g:M^m\to\R$. Since $\pi_P\left(f\left(M\right)\right)\subset B_P[R]$, 
we have that
$g^*\leq\psi_b\left(R\right)$. Let $\phi:M^m\to\left[0,+\infty\right)$ be given by
\beeqs
\phi\left(x\right)=\exp\left(\int_0^{\left|z\left(x\right)\right|}\frac{ds}{\varsigma\left(s\right)}\right),
\eneqs
where $\left|\,\right|$ stands for the distance function to the pole of $Q^l$. 
Since $f$ is proper and $\pi_P\left(f\left(M\right)\right)\subset B_P[R]$, 
then the function $\left|z\left(x\right)\right|$ satisfies 
$\left|z\left(x\right)\right|\to+\infty$ as $x\to\infty$. By hypothesis we 
have that $\int_0^{+\infty}1/\varsigma\left(s\right)ds=+\infty$ so that 
$\phi\left(x\right)\to+\infty$ as $x\to\infty$. 
We let $x_0\in M^m$ with $\pi_P\left(f\left(x_0\right)\right)\neq o$ and set
\beeqs
g_k\left(x\right)=\frac{g\left(x\right)-g\left(x_0\right)+1}{\phi\left(x\right)^{1/k}}.
\eneqs
Thus $g_k\left(x_0\right)>0$, and since $g^*\leq\psi_b\left(R\right)<+\infty$ and 
$\phi\left(x\right)\to+\infty$ as $x\to\infty$, we have that 
$\limsup_{x\to\infty}g_k\left(x\right)\leq0$. Hence $g_k$ attains a positive 
absolute maximum at a point $x_k\in M^m$. This procedure yields a 
sequence $\left\{x_k\right\}_{k\in\N}$ such that (passing to a subsequence
if necessary) $g\left(x_k\right)$ converges to $g^*$. First suppose that 
$x_k\to\infty$ as $k\to+\infty$. Since $g_k$ attains a maximum at $x_k$, 
we have $\grad^Mg_k\left(x_k\right)=0$ and 
$\hess^Mg_k\left(x_k\right)\left(X,X\right)\leq0$ for every 
$X\in T_{x_k}M$. This yields
\beeq
\grad^Mg\left(x_k\right)=\frac{g\left(x_k\right)-g\left(x_0\right)+1}
{k\phi\left(x_k\right)}\grad^M\phi\left(x_k\right)
\eneq
and
\beeq\label{hesg}
\begin{aligned}
\hess^Mg\left(x_k\right) & \leq\frac{g\left(x_k\right)-g\left(x_0\right)+1}{k\phi\left(x_k\right)}\left(\hess^M\phi\left(x_k\right)+\left(\frac{1}{k}-1\right)\frac{1}{\phi\left(x_k\right)}d\phi\otimes d\phi\right)\\
              & \leq\frac{g\left(x_k\right)-g\left(x_0\right)+1}{k\phi\left(x_k\right)}\hess^M\phi\left(x_k\right).
\end{aligned}
\eneq
Since $\phi\left(x\right)=\zeta\left(z\left(x\right)\right)$, where 
$\zeta\left(z\right)=\exp\left(\int_0^{\left|z\right|}ds/\varsigma\left(s\right)\right)\virg z\in Q^l$, from \eqref{hc} we have that
\beeq\label{heph}
\hess^M\phi(x)(X,X)=\hess^Q\zeta(z(x))(X_Q,X_Q)+
\langle\grad^Q\zeta(z(x)),\al_x(X,X)\rangle
\eneq
for all vectors $X\in T_xM$, where $X_Q=z_*X$. Also observe that
\beeqs
\grad^Q\zeta\left(z\right)=\frac{\zeta\left(z\right)}{\varsigma\left(\left|z\right|\right)}\grad^Q\left|z\right|,
\eneqs
and then
\beeq\label{gradphi}
\grad^M\phi\left(x\right)=\frac{\phi\left(x\right)}{\varsigma\left(\left|z\left(x\right)\right|\right)}\left(\grad^Q\left|z\left(x\right)\right|\right)^T.
\eneq
Thus, for every $X\in T_xM$ such that $X_Q=0$, it follows from \eqref{heph}
that
%\small
\begin{eqnarray*}
\hess^M\phi\left(x\right)\left(X,X\right) &=& 
\frac{\phi\left(x\right)}{\varsigma\left(\left|z\left(x\right)\right|\right)}
\langle\grad^Q\left|z\left(x\right)\right|,\al_x\left(X,X\right)\rangle \\ 
&\leq& \frac{\phi\left(x\right)}{\varsigma\left(\left|z\left(x\right)\right|\right)}
\left\|\al_x\left(X,X\right)\right\|.
\end{eqnarray*}
%\normalsize
Therefore, by \eqref{bsff} we obtain that
\beeq\label{hephe}
\frac{1}{\phi\left(x\right)}\hess^M\phi\left(x\right)\left(X,X\right)\leq\frac{\left\|\al_x\left(X,X\right)\right\|}{\varsigma\left(\left|z\left(x\right)\right|\right)}\leq\left\|X\right\|^2
\eneq
for every $X\in T_xM$ with $X_Q=0$.
Given $W_{x_k}\subset T_{x_k}M$ with $\dim W_{x_k}>p+l$, we have that the subspace $V_k=f_*^{-1}\left(f_*W_{x_k}\cap T_{y_k}P\right)$ has $\dim V_{x_k}\geq\dim W_{x_k}-l>p$ and $f_*\left(V_k\right)$ is orthogonal to $T_{z_k}Q$. Then, $X_Q=0$ for every $X\in V_{x_k}$, and from \eqref{hesg} and \eqref{hephe} 
we get that
\begin{eqnarray*}
\hess^Mg\left(x_k\right)\left(X,X\right) &\leq& 
\frac{g\left(x_k\right)-g\left(x_0\right)+1}{k\phi\left(x_k\right)}\hess^M\phi\left(x_k\right)\left(X,X\right) \\ 
&\leq& \frac{\psi_b\left(R\right)+1}{k}\left\|X\right\|^2,
\end{eqnarray*}
for every $X\in V_{x_k}$. Moreover, using Theorem \ref{hcth}, we also have here that
\beeq\label{hct2}
\hess^Mg\left(x\right)\left(X,X\right)\geq\psi_b'\left(r\left(y\left(x\right)\right)\right)\left(C_b\left(r\left(y\left(x\right)\right)\right)\left\|X\right\|^2-\left\|\al_x\left(X,X\right)\right\|\right)
\eneq
for every $X\in V_{x_k}$, since $X_P=X$. Therefore, we obtain that
\begin{eqnarray*}
\frac{\psi_b\left(R\right)+1}{k}\left\|X\right\|^2 &\geq& 
\hess^Mg\left(x_k\right)\left(X,X\right) \\ 
&\geq& 
\psi_b'\left(r_k\right)\left(C_b\left(r_k\right)\left\|X\right\|^2-\left\|\al_{x_k}
\left(X,X\right)\right\|\right)
\end{eqnarray*}
for every $x_k$ and every $X\in V_{x_k}$, where as usual $r_k=r\left(y_k\right)$. Hence
\beeqs
\left\|\al_{x_k}\left(X,X\right)\right\|\geq\left(C_b\left(r_k\right)-\frac{\psi_b\left(R\right)+1}{k\psi_b'\left(r_k\right)}\right)\left\|X\right\|^2
\eneqs
with
\beeqs
C_b\left(r_k\right)-\frac{\psi_b\left(R\right)+1}{k\psi_b'\left(r_k\right)}>0
\eneqs
for $k$ sufficiently large. Reasoning now as in the last part of the proof of Theorem \ref{math}, there
exists a plane $\sigma_k\subset V_{x_k}$ such that, by the Gauss equation
\beeqs
K_f\left(\sigma_k\right)=K_{\al}\left(\sigma_k\right)\geq\left(C_b\left(r_k\right)-\frac{\psi_b\left(R\right)+1}{k\psi_b'\left(r_k\right)}\right)^2,
\eneqs
and \eqref{oe} and \eqref{oie} follow by letting $k\to+\infty$ as in the last part 
of the proof of Theorem \ref{math}.
To finish the proof of Theorem \ref{sffcth}, we need to consider the case where the sequence $\left\{x_k\right\}_{k\in\N}\subset M^m$ remains in a compact set. In that case, passing to a subsequence if necessary, we may
assume that $x_k\to x_\infty\in M^m$ and $g$ attains its absolute maximum at $x_\infty$. Thus $\hess^Mg\left(x_\infty\right)\left(X,X\right)\leq0$ for all $X\in T_{x_\infty}M$. In particular, if follows from \eqref{hct2} that for every $X\in V_{x_\infty}$
\beeqs
0\geq\hess^Mg\left(x_\infty\right)\left(X,X\right)\geq\psi_b'\left(r_\infty\right)\left(C_b\left(r_\infty\right)\left\|X\right\|^2-\left\|\al_{x_\infty}\left(X,X\right)\right\|\right).
\eneqs
Therefore
\beeqs
\left\|\al_{x_\infty}\left(X,X\right)\right\|\geq C_b\left(r_\infty\right)\left\|X\right\|^2.
\eneqs
By applying Lemma \ref{olwv} to $\al|_{V_{x_\infty}\times V_{x_\infty}}:V_{x_\infty}\times V_{x_\infty}\to N_fM\left(x_\infty\right)$ and reasoning again as in the last part of the proof of Theorem \ref{math}, we have that there exists a plane $\sigma_\infty\subset V_{x_\infty}$ such that, by the Gauss equation,
\beeqs
K_f\left(\sigma_\infty\right)=K_{\al}\left(\sigma_\infty\right)\geq C_b^2\left(r_\infty\right),
\eneqs
and \eqref{oe} follows. Again \eqref{oie} follows as in the last part of the proof of Theorem \ref{math}.
\enpr

\bepr[Proof of Theorem \ref{ltfh}]
We proceed as in the proof of Theorem \ref{sffcth} to obtain a sequence $\left\{x_k\right\}_{k\in\N}$ such that $g\left(x_k\right)$ converges to $g^*$ and satisfying
\beeq\label{gradgh}
\grad^Mg\left(x_k\right)=\frac{g\left(x_k\right)-g\left(x_0\right)+1}{k\phi\left(x_k\right)}\grad^M\phi\left(x_k\right)
\eneq
and
\beeq\label{hessMh}
\hess^Mg\left(x_k\right)\leq\frac{g\left(x_k\right)-g\left(x_0\right)+1}{k\phi\left(x_k\right)}\hess^M\phi\left(x_k\right).
\eneq
Recall that (see \eqref{gradphi})
\beeq\label{gradphih}
\grad^M\phi\left(x\right)=\frac{\phi\left(x\right)}{\varsigma\left(\left|z\left(x\right)\right|\right)}\left(\grad^Q\left|z\left(x\right)\right|\right)^T.
\eneq
Let us first consider the case where $x_k\to\infty$ as $k\to+\infty$. From \eqref{gradgh} and \eqref{gradphih} we know that
\beeqs
\left\|\grad^Mg\left(x_k\right)\right\|\leq\frac{g^*+1}{k}\frac{1}{\varsigma\left(\left|z_k\right|\right)}\leq\frac{\psi_b\left(R\right)+1}{k}\frac{1}{\varsigma\left(\left|z_k\right|\right)}.
\eneqs
Since $f$ is proper and $\pi_P\left(f\left(M\right)\right)\subset B_P(R)$, 
then $\left|z_k\right|\to+\infty$ as $k\to+\infty$. Therefore, taking into account that $\limsup_{t\to+\infty}1/\varsigma\left(t\right)<+\infty$ we obtain from here that
\beeq
\lim_{k\to+\infty}\left\|\grad^Mg\left(x_k\right)\right\|=0.
\eneq
Observe that
\beeqs
\grad^Nh\left(f\left(x\right)\right)=\psi_b'\left(r\left(y\right)\right)\grad^Pr\left(y\right)=\grad^Mg\left(x\right)+\left(\grad^Nh\left(f\left(x\right)\right)\right)^\perp,
\eneqs
where $y=y\left(x\right)$. Therefore,
\beeq\label{leq1}
\psi_b'\left(r_k\right)^2=\left\|\grad^Mg\left(x_k\right)\right\|^2+\left\|\left(\grad^Nh\left(f\left(x_k\right)\right)\right)^\perp\right\|^2,
\eneq
and making $k\to+\infty$ here we obtain that
\beeqs
\lim_{k\to+\infty}\left\|\left(\grad^Nh\left(f\left(x_k\right)\right)\right)^\perp\right\|=\psi_b'\left(r^*\right)>0,
\eneqs
which implies that
\beeqs
\left(\grad^Nh\left(f\left(x_k\right)\right)\right)^\perp\neq0
\eneqs
for $k$ sufficiently large. As in the proof of Theorem \ref{sffcth}, since 
$n-l>1$, given $W_{x_k}\subset T_{x_k}M$ with $\dim W_{x_k}>l+1$, we have that $V_k=f_*^{-1}\left(f_*W_{x_k}\cap T_{y_k}P\right)$ has $\dim V_{x_k}\geq\dim W_{x_k}-l>1$ and $f_*\left(V_k\right)$ is orthogonal to $T_{z_k}Q$. Then, using Theorem \ref{hcth}, we also have that
\beeq\label{hcteq}
\hess^Mg\left(x_k\right)\left(X,X\right)\geq\psi_b'\left(r_k\right)\left(C_b\left(r_k\right)\left\|X\right\|^2-\left\|\al_{x_k}\left(X,X\right)\right\|\right)
\eneq
for every $X\in V_{x_k}$, since $y_*X=X$. On the other hand, we also know from \eqref{hessMh} that
\beeq\label{leq3}
\begin{aligned}
\hess^Mg\left(x_k\right)\left(X,X\right) & \leq\frac{\psi_b\left(R\right)+1}{k}\frac{\hess^M\phi\left(x_k\right)\left(X,X\right)}{\phi\left(x_k\right)}\\
                   & =\frac{\psi_b\left(R\right)+1}{k}\frac{1}{\varsigma\left(\left|z_k\right|\right)}\langle\grad^Q\left|z_k\right|,\al_{x_k}\left(X,X\right)\rangle
\end{aligned}
\eneq
for every $X\in T_{x_k}M$. Since we are in codimension one and $\left(\grad^Nh\left(f\left(x_k\right)\right)\right)^\perp\neq0$ (for $k$ large enough), then
\beeq\label{leq2}
\al_{x_k}\left(X,X\right)=\lambda_k\left(X,X\right)\left(\grad^Nh\left(f\left(x_k\right)\right)\right)^\perp
\eneq
for a real function $\lambda_k$. Now observe that
\beeqs
\begin{aligned}
\langle\grad^Q\left|z_k\right|,\al_{x_k}\left(X,X\right)\rangle&=\lambda_k\left(X,X\right)\langle\grad^Q\left|z_k\right|,\left(\grad^Nh\left(f\left(x_k\right)\right)\right)^\perp\rangle\\
&=\lambda_k\left(X,X\right)\langle\grad^Q\left|z_k\right|,\grad^Mg\left(x_k\right)\rangle
\end{aligned}
\eneqs
because of $\langle\grad^Q\left|z_k\right|,\grad^Pr\left(y_k\right)\rangle=0$. Therefore,
\beeqs
\begin{aligned}
\langle\grad^Q\left|z_k\right|,\al_{x_k}\left(X,X\right)\rangle&\leq\left|\lambda_k\left(X,X\right)\right|\left\|\grad^Mg\left(x_k\right)\right\|\\
&\leq\left|\lambda_k\left(X,X\right)\right|\frac{\psi_b\left(R\right)+1}{k}\frac{1}{\varsigma\left(\left|z_k\right|\right)}.
\end{aligned}
\eneqs
On the other hand, from our hypothesis \eqref{hyphyp} we know that
\beeqs
\left\|\al_x\left(X,X\right)\right\|\leq\varsigma^2\left(\left|z\left(x\right)\right|\right)\left\|X\right\|^2,
\eneqs
and from \eqref{leq1} and \eqref{leq2} we have that
\beeqs
\left\|\al_{x_k}\left(X,X\right)\right\|=\left|\lambda_k\left(X,X\right)\right|\sqrt{\psi_b'\left(r_k\right)^2-\left\|\grad^Mg\left(x_k\right)\right\|^2}\leq\varsigma^2\left(\left|z_k\right|\right)\left\|X\right\|^2.
\eneqs
That is,
\beeqs
\frac{\left|\lambda_k\left(X,X\right)\right|}{\varsigma\left(z_k\right)}\leq\frac{\varsigma\left(z_k\right)\left\|X\right\|^2}{\sqrt{\psi_b'\left(r_k\right)^2-\left\|\grad^Mg\left(x_k\right)\right\|^2}}.
\eneqs
It follows from here that
\beeqs
\langle\grad^Q\left|z_k\right|,\al_{x_k}\left(X,X\right)\rangle\leq\frac{\psi_b\left(R\right)+1}{k}\frac{\varsigma\left(z_k\right)\left\|X\right\|^2}{\sqrt{\psi_b'\left(r_k\right)^2-\left\|\grad^Mg\left(x_k\right)\right\|^2}}
\eneqs
for every $X\in T_{x_k}M$, so that by \eqref{leq3} we get
\beeq\label{leq4}
\hess^Mg\left(x_k\right)\left(X,X\right)\leq\left(\frac{\psi_b\left(R\right)+1}{k}\right)^2\frac{\left\|X\right\|^2}{\sqrt{\psi_b'\left(r_k\right)^2-\left\|\grad^Mg\left(x_k\right)\right\|^2}}.
\eneq
Therefore, from \eqref{hcteq} and \eqref{leq4} we have that
\begin{eqnarray*}
\psi_b'\left(r_k\right)\left(C_b\left(r_k\right)\left\|X\right\|^2-\left\|\al_{x_k}\left(X,X\right)\right\|\right)\leq\left(\frac{\psi_b\left(R\right)+1}{k}\right)^2\frac{\left\|X\right\|^2}{\sqrt{\psi_b'\left(r_k\right)^2-\left\|\grad^Mg\left(x_k\right)\right\|^2}}
\end{eqnarray*}
for every $X\in V_{x_k}$. Hence
\beeqs
\left\|\al_{x_k}\left(X,X\right)\right\|\geq\left(C_b\left(r_k\right)-\frac{\left(\psi_b\left(R\right)+1\right)^2}{k^2\psi_b'\left(r_k\right)\sqrt{\psi_b'\left(r_k\right)^2-\left\|\grad^Mg\left(x_k\right)\right\|^2}}\right)\left\|X\right\|^2,
\eneqs
with
\small
\beeqs
\lim_{k\to+\infty}\left(C_b\left(r_k\right)-\frac{\left(\psi_b\left(R\right)+1\right)^2}{k^2\psi_b'\left(r_k\right)\sqrt{\psi_b'\left(r_k\right)^2-\left\|\grad^Mg\left(x_k\right)\right\|^2}}\right)=C_b\left(r^*\right)\geq C_b\left(R\right)>0.
\eneqs
\normalsize
Reasoning now as in the last part of the proof of Theorem \ref{math}, there exists a plane $\sigma_k\subset V_{x_k}$ such that, by the Gauss equation,
\beeqs
\begin{aligned}
K_f\left(\sigma_k\right)\geq\left(C_b\left(r_k\right)-\frac{\left(\psi_b\left(R\right)+1\right)^2}{k^2\psi_b'\left(r_k\right)\sqrt{\psi_b'\left(r_k\right)^2-\left\|\grad^Mg\left(x_k\right)\right\|^2}}\right)^2,
\end{aligned}
\eneqs
and \eqref{oe} and \eqref{oie} follow by letting $k\to+\infty$ as in the last part of the proof of Theorem \ref{math}.
Finally, in the case where the sequence $\left\{x_k\right\}_{k\in\N}\subset M^n$ remains in a compact subset of $M^n$, and passing to a subsequence if necessary, we may assume that $x_k\to x_\infty\in M^n$ and $g$ attains its absolute maximum at $x_\infty$. Thus, $\hess^Mg\left(x_\infty\right)\left(X,X\right)\leq0$ for all $X\in T_{x_\infty}M$. Therefore, it follows again from Theorem \ref{hcth} that for every $X\in V_{x_\infty}$,
\beeqs
0\geq\hess^Mg\left(x_\infty\right)\left(X,X\right)\geq\psi_b'\left(r_\infty\right)\left(C_b\left(r_\infty\right)\left\|X\right\|^2-\left\|\al_{x_\infty}\left(X,X\right)\right\|\right).
\eneqs
The proof now finishes as in Theorem \ref{sffcth}.
\enpr
\section{Notes}
The idea of the proof that any compact surface in $\R^3$ must have a point of positive
Gauss curvature was first taken up by Tompkins \cite{MR1546106}, who showed that there is no
isometric immersion $f:M^n\to\R^{2n-1}$ if $M^n$ is compact and flat. This result inspired
the seminal paper of Chern-Kuiper \cite{MR0050962}, where Lemma \ref{olwv} was proved for dimensions
$n=2,\,3$ and conjectured to be true for any dimension. This conjecture was proved by
Otsuki \cite{MR0067550} for $\lambda=0$ who, consequently, obtained Theorem \ref{tols} for all dimensions.

The Chern and Kuiper result gave rise to a long series of works, among others,
by O'Neill \cite{MR0117683}, Stiel \cite{MR0212743}, Moore \cite{MR0377776}, Jorge-Koutroufiotis \cite{MR623135}, Pigola-Rigoli-Setti \cite{MR2116555} and, finally, by Al\'{i}as-Bessa-Montenegro \cite{MR2901222} who obtained Theorem \ref{abmt} on cylindrically bounded submanifolds.

The maximum principles used throughout this paper, as well as their related results, namely, Theorems \ref{oympct}, \ref{oywmpct} and \ref{oympsc}, are due to Pigola-Rigoli-Setti
\cite{MR2116555}. On the other hand, it was shown in \cite{MR3148611} that conditions (\emph{b}) and (\emph{c}) in Theorem \ref{oympct} can be replaced by the following equivalent although apparently stronger requirements:
\begin{itemize}
\item[(b)] $\left\|\grad^M\gamma\right\|\leq c$ for a constant $c>0$ outside a compact subset of $M^n$,
\item[(c)] $\hess^M\gamma\leq d\langle\,,\,\rangle$ for a constant $d>0$ outside a compact subset of $M^n$.
\end{itemize}
A similar observation holds for the Omori-Yau maximum principle.

Regarding complete hypersurfaces of nonpositive Ricci curvature, Leung \cite{MR707970} pioneered their study by proving Theorem \ref{veer} in the case $b=0$ and conjecturing that the assumption on the sectional curvature could be dropped. This, however, turns out not to be true, as shown by Nadirashvili's \cite{MR1419004} celebrated counterexample to both Hadamard's and Calabi-Yau's conjectures on negatively curved and minimal surfaces. After Leung's work, Smith \cite{MR729766} gave an answer for the case $b<0$ but with a non-sharp estimate (for having made the Hessian comparison to $\R^{n+1}$ instead of $\mathbb{H}_b^{n+1}$), and finally Veeravalli \cite{MR1775899} obtained Theorem \ref{veer}. 

\vspace{.2cm}

It is a natural question to ask whether Theorem \ref{svt} is still true in the limiting case, that is, when $R=\inj_P\left(o\right)=\frac{\pi}{2\sqrt{b}}$, where $\frac{\pi}{2\sqrt{b}}$ is replaced by $+\infty$ if $b\leq0$. This motivates the following conjecture.

\beco\label{ourconj}
Let $f:M^m\to N^{n+l}=P^n\times Q^l$ be an isometric immersion with 
codimension $p=n+l-m<m-l$ of a complete Riemannian manifold. Assume
that $R=\inj_P\left(o\right)=\frac{\pi}{2\sqrt{b}}$, where $\frac{\pi}{2\sqrt{b}}$ is replaced by $+\infty$ if $b\leq0$. If $K^{\rad}_P\leq b$ in $B_P[R]$, then
\beeq
\sup_M\min\left\{\max_{\sigma\subset W}K_f\left(\sigma\right):\dim W>p+l\right\}
\geq\max\left\{-b,0\right\}.
\eneq
Moreover,
\beeq
\sup_M\min\left\{\max_{\sigma\subset W}K_M\left(\sigma\right):\dim W>p+l\right\}\geq\max\left\{-b,0\right\}+\inf_{B_P[R]}K_P.
\eneq
\enco

It is not clear the extent to which the above conjecture is true, but an affirmative answer at least in the most classic cases, such as $P^n=\R^n$ and $l=0$, would have deep implications in the field of submanifolds with nonpositive extrinsic curvature. Indeed, Conjecture \ref{ourconj} in this case implies when $p=m-1$
that a complete Riemannian manifold with sectional curvature $K\leq-c<0$ cannot be immersed isometrically in $\R^{2m-1}$, a kind of Efimov's theorem in $n$ dimensions. In particular, this would give us the $m$-dimensional version of the classical theorem of Hilbert that the hyperbolic plane cannot be realized isometrically in $\R^3$.

There is yet another attempt to extend Efimov's theorem to higher dimensions in a different direction proposed independently by Reilly \cite{MR0474149} and Yau \cite{MR0370443} (see also \cite{MR645728} and Gromov \cite{MR864505}):

``\emph{There are no complete hypersurfaces in $\R^{n+1}$ with Ricci curvature $\leq-c$}'' and proved to be very true for $n=3$ and essentially true for $n>3$ by Smyth-Xavier \cite{MR914845}. Their main result seems to be inaccessible to techniques using the Omori-Yau maximum principle, and its proof relies on a purely geometric result on the principal curvatures of complete submanifolds of Euclidean space. Still in the case $P^n=\R^n$ and $l=0$, Conjecture \ref{ourconj} for $p=1$ would not only settle the above question at all, but also, in the same spirit of Corollary \ref{clvsc}, weaken the assumption that $\Ric\leq-c$ to $s\leq-c$.

%\bibliographystyle{plain}
%\bibliography{ourbib}
\end{document}